\documentclass[12pt,reqno]{amsart}
\usepackage{amsmath}
\usepackage{amssymb}
\usepackage{amsfonts}
\numberwithin{equation}{section}
\newcommand{\CC}{\mathbb{C}\,}

\newcommand{\Fp}{{\mathcal F_\varphi}}

\newcommand{\kl}{{\Bbbk_\lambda}}

\newcommand{\kol}{{\bf k}_{\lambda}}

\newcommand{\elll}{{\ell}}

\newcommand{\Hol}{\mathop{\rm Hol\,}}

\DeclareMathOperator{\dist}{dist}
\DeclareMathOperator{\card}{Card}
\DeclareMathOperator{\const}{Const}

\newtheorem{th1}{{\bf Theorem}}[section]

\newtheorem{thm}[th1]{{\bf Theorem}}
\newtheorem{lem}[th1]{{\bf Lemma}}

\theoremstyle{definition}

\begin{document}
\author[A. Borichev, Yu. Lyubarskii]
{Alexander Borichev, Yurii Lyubarskii}

\title[Riesz bases of reproducing kernels in Fock type
spaces]
{Riesz bases of reproducing kernels in Fock type
spaces}

\begin{abstract}
In a scale of Fock spaces $\Fp$ with radial weights $\varphi$ we study the existence
of Riesz bases of (normalized) reproducing kernels. We prove that these spaces possess such bases if and only if $\varphi(x)$ grows at most like $(\log x)^2$.
\end{abstract}

\thanks{A.B.  was partially supported by the ANR project DYNOP;
 Yu.L. was partially supported by the
Research Council of Norway grants 10323200 and  160192/V30.
}

\maketitle

\section{Introduction.}

Given an increasing   function $\varphi$ defined on $[0,+\infty)$,
we extend it to $\CC$ by $\varphi(z)=\varphi(|z|)$, and consider
the Fock type space 
$$
\Fp=\{f\in\Hol(\CC):\|f\|^2_\varphi= \int_{\CC}|f(z)|^2e^{-2\varphi(z)}dm(z)<\infty\},
$$
where $dm$ is the area Lebesgue measure. 

The Hilbert space $\Fp$ possesses  the bounded point evaluation property,
i.e.,  for each $\lambda \in \CC$,  the mapping $L_\lambda: f \mapsto f(\lambda)$  is a bounded linear functional in $\Fp$.  Therefore there exists
 ${\bf k}_\lambda={\bf k}^\varphi_\lambda\in \Fp$, the
reproducing kernel at $\lambda$ in $\Fp$:
\[
f(\lambda)=\langle f, {\bf k}_\lambda \rangle_\varphi, \qquad 
\lambda\in \CC,\ f\in\Fp,
\]
and we have
\[
\|L_\lambda\|_{\Fp\to \CC}= \| {\bf k}_\lambda\|_\varphi=
\big ({\bf k}_\lambda(\lambda)\big )^{1/2}. 
\] 

Let $\kl={\bf k}_\lambda/\|{\bf k}_\lambda\|_{\varphi}$
be the normalized reproducing kernel at $\lambda$. 
Given a sequence 
$\Lambda \subset\mathbb C$,
we say that $\{\Bbbk_{\lambda}\}_{\lambda \in \Lambda}$ is 
a Riesz basis in $\Fp$ if it is complete and for some $c,C>0$ we have
$$
c\sum_{\lambda \in \Lambda}|a_\lambda|^2\le \Bigl\|\sum_{\lambda \in \Lambda}
a_\lambda \Bbbk_{\lambda}\Bigr\|^2_{\varphi}
\le C\sum_{\lambda \in \Lambda}|a_\lambda|^2,  
$$
for each finite sequence $\{a_\lambda\} \subset \mathbb C$.
Equivalently $\{\Bbbk_{\lambda}\}_{\lambda\in \Lambda}$
is a linear isomorphic image
of an orthonormal basis in a separable Hilbert space.

In this article we study the following question: {\em for which $\varphi$ does the space 
$\Fp$  admit a Riesz basis of normalized  reproducing kernels?}

The above question can be reformulated in the 
classical terms of  interpolation
in Hilbert spaces of entire functions. 
Let  $X$ be such a space and 
let ${\bf  k}^X_\lambda$ stand for  the reproducing 
kernel at $\lambda$ 
in $X$. We say that a sequence $\Lambda \subset \CC$  is a 
{\em complete interpolating sequence for $X$} if, for each 
$\{a_\lambda\}\in \elll^2(\Lambda) $,  the interpolation problem
\begin{equation}
\label{eq03a}
\frac {f(\lambda)}{ \| {\bf k}^X_\lambda\|_X
}=a_\lambda, \qquad 
\lambda \in \Lambda,
\end{equation}
has a unique solution in $X$.

Standard duality arguments show that the system 
$\{ \kl \}_{\lambda \in \Lambda}$ is a Riesz basis  
in 
$X$  if and only if $\Lambda $  is a complete interpolating sequence for this space, %
$\kl = {\bf k}^X_\lambda/\| {\bf k}^X_\lambda\|_X$.

A canonical example here is the Paley-Wiener space, i.e. the space of Fourier transforms
of all functions from $L^2(-\pi,\pi)$. The set 
$\Lambda=\mathbb Z$ is a complete interpolating sequence for this space.
Furthermore, a notoriously difficult open problem by Nikolski is whether every model space 
$K_\Theta=H^2\ominus \Theta H^2$ ($\Theta$ is an inner function) has a Riesz basis of reproducing
kernels.
We refer also to the papers \cite{Lu,Is,Ly} studying expansions in exponential series on different domains of the complex plane. 
These problems can be reformulated in terms of bases
of reproducing kernels in Fock type spaces. In particular, the results of
\cite{Is} show that $\Fp$ with
$\varphi(x)=x-\frac 32\log x$ has no Riesz basis of (normalized) reproducing
kernels.

Starting from the results of Seip \cite{SW1}, 
it is known that the classical Fock space, $\Fp$ with
$\varphi(x)=x^2$, unlike the Paley--Wiener spaces, has no Riesz basis of (normalized) reproducing
kernels.

For more rapidly growing  
$\varphi$
the absence of such bases is established in \cite{BDK} (see Theorems~2.2 and 2.4 therein) under some natural regularity conditions on $\phi$.   On the other hand, if $\varphi(x)= \const\cdot \log |x| $, then    $\Fp$  becomes
a finite dimensional space of polynomials so that each  
$\Lambda \subset \CC$ with $\card (\Lambda)= \dim(\Fp) $ is obviously  a complete interpolating sequence for this space.

In this paper we assume that $\varphi(z)=\varphi(|z|)$ is a subharmonic function such that
\[
\log |x| = o(\varphi (x)), \qquad x\to \infty,
\]
hence $\dim \Fp= \infty$, and 
study what happens
if $\varphi(x)$ grows less rapidly than $x^2$.

It turns out that for $\varphi(x)=(\log x)^\beta$, $1<\beta\le 2$, the spaces $\Fp$ 
still have Riesz bases of normalized reproducing kernels. 
On the other hand, for   $\varphi$ such that
$(\log x)^2=o(\varphi(x))$, $x\to\infty$, the spaces $\Fp$ have no
such bases (again assuming that $\phi$ satisfies some natural regularity conditions). 
Roughly speaking, the reason for this is that the local scale function 
$\rho(z)=(\Delta\varphi(z))^{-1/2}$ is $o(|z|)$, $z\to\infty$, if $(\log x)^2=o(\varphi(x))$,
$x\to\infty$, and is comparable to $z$, if $\varphi(x)=(\log x)^2$.

Precise formulations of our results are given in Section 2.
We study the case when both $(\log x)^2=o(\varphi(x))$ and $\varphi(x)=O(x^2)$
in Section 3.
There we apply a theorem from \cite{LM}  in order to
approximate   $\exp(\varphi)$
by the modulus of an entire function using  discretization of the
Riesz measure $d\mu_\varphi =\triangle \varphi(z)\,dm(z)$ 
of the subharmonic function $\varphi$,
  and then use an argument by Seip 
from \cite[Lemma 6.2]{SW1}. In Section 4 we deal with the borderline
case $\varphi(x)=(\log x)^2$, and again use  approximation 
of $\exp(\varphi)$ by the modulus of an entire function.
Finally, our argument in Section 5 dealing with the case
$\varphi(x)=x^\beta$, $1<\beta< 2$, is essentially a real variable
one using Legendre transform estimates.

In Section 2 we obtain asymptotic estimates (Lemma~\ref{le1a}, Lemma~\ref{le4}) 
of the norm of the reproducing kernel coinciding with those by Holland--Rochberg 
in \cite{HOR}. We cannot just refer to \cite{HOR} since the conditions  
imposed on $\varphi$ there are not satisfied in our situation.
\medskip

\noindent{\bf Acknowledgment}  This work  began when the second author visited
the University  of Aix-Marseille. He thanks the University for hospitality and support.  

The authors are thankful to the referees for numerous comments.

\section{Formulation of the  results}
\label{sect2}

\subsection*{Part A} In this subsection we deal with regular functions $\varphi(x)$ 
growing more rapidly than $(\log x)^2$, but less rapidly than 
$x^2$. We assume that $\varphi(z)=\varphi(|z|)$
is $C^2$-smooth and subharmonic on $\CC$, and set 
$$
\rho(z)=\bigl(\triangle \varphi(z)\bigr)^{-1/2}=
\Bigl(\frac{\varphi'(r)}{r}+\varphi''(r)\Bigr)^{-1/2},
$$
where $r=|z|$.  

The function $\rho(z)$  defines a natural scale with respect to the Riesz
measure $\mu_\varphi$ of   $\varphi$, i.e. 
$\mu_\varphi(\{\zeta: |\zeta-z|\leq \rho(z)\})\asymp 1$  for all $z\in \CC$.
Here and in what follows, the notation $A(s)\prec B(s)$ for $s$ in some set $S$ 
means that the ratio $A(s)/B(s)$ of the two positive functions $A(s)$ and $B(s)$ 
is bounded from above by a  positive constant independent of $s$ in $S$. 
We write $A(s) \succ B(s)$ if $B(s)\prec A(s)$ and $A(s)\asymp B(s)$  
if both $A(s)\prec B(s)$  and $B(s)\prec A(s)$.

In the borderline cases $\varphi(x)=x^2$  and $\varphi(x)=(\log x)^2$  we have respectively
$\rho(r)=\const$  and $\rho(r)= \const \cdot r $.

In this subsection we assume that
\begin{equation}
0< \inf_{r>0} \rho(r), \ \text{and}  \     \rho(r)=o(r), \quad r\to\infty,\\\label{01}
\end{equation}
and also 
\begin{equation}
\left.
\begin{gathered}
 \rho(r+ \rho(r))=(1+o(1))\rho(r),\quad  r\to\infty,\\
\rho(2r)\asymp \rho(r),\qquad r>0. 
\end{gathered}
\right\}
\label{01a}
\end{equation}

\begin{lem} Given $w\in\mathbb C$, there exists a function $\Phi_w$ analytic
in the disc $D_w=\{z\in\mathbb C:|z-w|<\rho(w)\}$ and such that
$$
|\Phi_w(z)|\asymp  e^{\varphi(z)},\qquad z\in D_w.
$$
\label{le0}
\end{lem}

\begin{lem} There exists an entire function $F$ such that
\begin{equation}
\label{eq01}
|F(z)|\asymp e^{\varphi(z)}\cdot \frac{\dist(z,W)}{\rho(z)},\qquad
z\in \CC,
\end{equation}
where $W$ is the zero set of $F$,
and 
\begin{gather*}
 \dist(w,W\setminus\{w\})\succ \rho(w), \quad w \in W, \\
\dist(z,W)\prec \rho(z), \quad z\in\CC. 
\end{gather*}
\label{le1}
\end{lem}

\begin{lem} 
$$
\|{\bf k }_z\|^2_{\varphi}\asymp e^{2\varphi(z)}/\rho^2(z), 
\qquad z \in \CC.
$$
\label{le1a}
\end{lem}

\begin{lem} Let a sequence  $\Lambda\subset \CC$ be such that 
$\{\kl \}_{\lambda \in \Lambda}$  is a Riesz basis in $\Fp$.
Then   
\begin{itemize}
\item[(a)] $\dist(\lambda,\Lambda\setminus
\{\lambda\})\succ \rho(\lambda)$, $\lambda \in \Lambda$,
\item[(b)] $\dist(z,\Lambda)\prec \rho(z)$, $z\in\CC$.
\end{itemize}
\label{le2}
\end{lem}

\begin{thm} Under conditions \eqref{01}, \eqref{01a} the space $\Fp$
has no Riesz bases of normalized reproducing kernels.
\label{te1}
\end{thm}
\medskip

\subsection*{Part B} In this subsection we assume that $\varphi(r)=(\log^+r)^2$.

\begin{lem} Let $\Lambda = \{\exp(\frac{n+1}2 + i\theta_n)\}_{n\ge 0}$,
where $\theta_n$ are arbitrary real numbers.  The product 
$$
E(z)=\prod_{\lambda \in \Lambda}\Bigl(1-\frac z{\lambda}\Bigr)
$$
converges uniformly on compact sets in $\CC$ and satisfies the estimate
\begin{equation}
\label{eq11}
|E(z)|\asymp e^{\varphi(z)}\cdot \frac{\dist(z,\Lambda)}{|z|^{3/2}},\qquad
z\in \CC.
\end{equation}
\label{le3}
\end{lem}

\begin{lem} 
$$
\|{\bf k}_z \|^2_{\varphi}\asymp e^{2\varphi(z)}/(1+|z|^2), 
\qquad z\in\mathbb C.
$$
\label{le4}
\end{lem}

\begin{thm} Let $\varphi(r)=(\log^+r)^2$, and let $\Lambda$
be   as in Lemma~{\rm \ref{le3}}. Then 
$\{\kl\}_{\lambda \in \Lambda}$ is a Riesz basis in $\Fp$.
\label{te2}
\end{thm}
\medskip

\subsection*{Part C}  In this subsection we consider the case 
$\varphi(r)=(\log^+r)^{1+\delta}$  for $0<\delta<1$.
\smallskip

 Denote $w_n=\log\|z^n\|^2_{\varphi}$, $n\ge 0$.
Then
$$
\kol(z)=\sum_{n\ge 0}\bar \lambda^nz^ne^{-w_n} 
$$
and
$$
\|\kol\|^2_{\varphi}=\sum_{n\ge 0}|\lambda|^{2n}e^{-w_n}.
$$

\begin{lem} For some $c>0$ we have
$$
 w_n= c(n+1)^{1+1/\delta} +  O\bigl(\log n\bigr),\qquad n> 0.
$$
\label{le5}
\end{lem}

Denote $r_0=0$, $r_n=\exp[(w_{n+1}-w_{n-1})/4]$, $n\ge 1$.

\begin{thm} Let $0<\delta<1$, $\varphi(r)=(\log^+r)^{1+\delta}$, 
and let $\lambda_n=r_ne^{i\theta_n}$ with arbitrary 
real $\theta_n$.
 Then $\{\Bbbk_{\lambda_n}\}_{n\ge 0}$ is a Riesz basis in $\Fp$.
\label{te3}
\end{thm}
\bigskip

\section{Proofs. Part A.}
\label{sect3}

\begin{proof}[Proof of Lemma~\rm\ref{le0}] (See also \cite[Lemma 4.1]{BDK}.)
By \eqref{01a} we know that $\rho(z)\asymp \rho(w)$, $z\in D_w $.
Set $H(\zeta)=\varphi(w+\zeta \rho(w))$, $|\zeta|\le 1$. Then
$\Delta H(\zeta)\asymp 1$, $|\zeta|\le 1$.
Next we define
$$
G(z)=\int_{\mathbb D}\log \Bigl|\frac{\zeta-z}{1-\overline{\zeta} z}\Bigr|\,\Delta H(\zeta)\,
dm(z),\qquad |z|\le 1.
$$
Then $|G(z)|\prec 1$, $|z|\le 1$  uniformly with respect $w\in \CC$, 
and $H_1=H-G$ is real and harmonic.
Denote by $\widetilde H_1$ the harmonic conjugate to $H_1$, $\widetilde H_1(0)=0$,
set $H_0=H_1+i\widetilde H_1$, and define 
$$
\Phi_w(z)=\exp H_0\bigl((z-w)/\rho(w)\bigr).
$$
Since $\log|\Phi_w(z)|-\varphi(z)=-G\bigl((z-w)/\rho(w)\bigr)$, $z\in D_w$, the proof is completed.
\end{proof}

\begin{proof}[Proof of Lemma~\rm\ref{le1}] 
This lemma is a special case of Theorem 3 in \cite{LM}. 
We  describe  just the idea of its  proof. It    relies 
on the   atomization procedure for the measure $\mu_\varphi$.  In this procedure
the complex plane  is  decomposed into a disjoint union of pieces $\omega$ of 
$\mu_\varphi$-measure two each, and then an atomized measure 
$\mu^{(a)}_\varphi$  is constructed; this measure is the sum of  discrete unit masses,
two masses are situated in each piece $\omega$ in such a way that the first 
two moments of their sum coincide with the corresponding moments
of $\mu_\varphi | \omega$. We refer the reader to \cite{Yu,Or}
for other implementation of atomization techniques.
\end{proof}

\begin{proof}[Proof of Lemma~\rm\ref{le1a}] 
We use the fact that 
$\|{\bf k}_z\|_\varphi=\|L_z\|_{\Fp \to \CC}$. 
Given $z \in\CC$, take $w,w'\in W$, $w\ne w'$, such that
$|z-w|=\dist(\lambda,W)$, $|z-w'|\asymp \rho(z)$,
and consider the function $G=F/[(\cdot-w)(\cdot-w')]$,
where $F$ is the function from Lemma~\ref{le1}, and $W$ is its zero set.

Then
\begin{gather}
|G(z)|\asymp e^{\varphi(z)}/\rho^2(z),\notag\\
\|G\|^2_{\varphi}\prec \int_{\CC}\frac{dm(\zeta)}{\rho^4(z)+|\zeta-z|^4}
\asymp \frac 1{\rho^2(z)}.\label{03}
\end{gather}
Therefore,
$$
\|{\bf k}_z^\varphi \|^2_\varphi \succ  e^{2\varphi(z)}/\rho^2(z).
$$
\smallskip

Now take any $f\in \Fp$, $\|f\|_{\varphi}=1$, and define $\Phi_z$, $D_z$
as in  Lemma~\ref{le0}. We have
\[
\int_{D_z} \bigl | f(\zeta)/ \Phi_z(\zeta) \bigr |^2 dm(\zeta) 
\prec \|f\|^2_\varphi = 1,
\]
 and the mean value theorem yields 
 \[
 |f(z)| \prec  |\Phi_z(z)|/\rho(z) \asymp 
 e^{\varphi(z)}/\rho(z). 
 \]
Thus,
$$
\|{\bf k}_z\|_{_\varphi}\asymp e^{\varphi(z)}/\rho(z).
$$
\end{proof}
\medskip

Now let, for a sequence $\Lambda  \subset \CC $,  the system 
$\{{\Bbbk}_\lambda \}_{\lambda \in \Lambda}$  be  a Riesz basis
in $\Fp$.
Then the mapping 
\[
f\mapsto \{ \langle f, {\bf k}_\lambda \rangle_{\varphi} \}_{\lambda\in \Lambda} =
    \{f(\lambda)\}_{\lambda\in \Lambda}
\]
is an isomorphism between $\Fp$  and the space
 $$
 \elll^2(1/\|{\bf k}_\lambda\| )=\Bigl\{\{c_\lambda\}_{\lambda \in \Lambda}:
 \ \|\{c_\lambda\}\|_ {\elll^2(1/\|{\bf k}_\lambda\| )}=
\sum_{\lambda \in \Lambda}|c_\lambda|^2/\|{\bf k}_\lambda\|^2_\varphi  < \infty \Bigr\}.
$$   
In particular, the interpolating problem 
\begin{equation}
\label{eq17}
f(\lambda)=c_\lambda, \qquad \lambda\in \Lambda, \ f\in \Fp
\end{equation}
has a unique solution for each $\{c_\lambda\}\in \elll^2(1/\|{\bf k}_\lambda\| )$ and 
\begin{equation}
\label{eq05}
\|f\|_\varphi \asymp \|f|\Lambda \|_ {\elll^2(1/\|{\bf k}_\lambda\| )}, \qquad f\in \Fp.
\end{equation}

Fix $\lambda  \in \Lambda$ and consider the function $f_\lambda \in \Fp$, solving the 
interpolating problem $f_\lambda(\mu)= \delta_{\lambda, \mu}$, $\mu\in \Lambda$,
here $\delta_{\lambda, \mu}$  is the Kronecker delta function.  Then the function 
$E(z)=(z-\lambda)f_\lambda(z)$ vanishes precisely on $\Lambda$ 
(otherwise, $\Lambda$ would not be a uniqueness set), and the solution to the problem  \eqref{eq17}  has the form 
\[
f(z) = \sum_{\lambda\in \Lambda}   c_\lambda 
      \frac{E(z)}{E'(\lambda)(z-\lambda)},
\]
the sum being convergent in $\Fp$.
 \medskip

\begin{proof}[Proof of Lemma~\rm\ref{le2}] 

(a) Suppose that $0<|\lambda -\lambda' |\le \rho(\lambda')/N$, 
$\lambda,\lambda'\in \Lambda$.
Denote $D=D_{\lambda'}$, and for $\zeta \in \CC$ set 
$$
f(z)=\frac{E(z)}{z-\lambda},\qquad
f_1(z)=\frac{E(z)}{(z-\lambda)(z-\lambda')}(z-\lambda'-\zeta \rho(\lambda')/2).
$$
Then for some $\zeta\in \CC$, $|\zeta|=1$, we have
\begin{equation}
\|f_1\|_\varphi \asymp \|f\|_\varphi.
\label{02}
\end{equation}
Indeed, the corresponding integrals
outside  $D$ are equivalent. Furthermore, let
$$
g(z)=\frac{f(z)}{z-\lambda'}\cdot \frac 1{\Phi_{\lambda'}(z)}.
$$
Then
$$
\int_D|f(z)|^2 e^{-2\varphi(z)}dm(z) \asymp
      \int_D|g(z)|^2|z-\lambda'|^2dm(z),
$$
and, for an appropriate $\zeta$, $|\zeta|=1$,
\begin{multline*}
\int_D|f_1(z)|^2 e^{-2\varphi(z)}dm(z) \asymp
      \int_D|g(z)|^2|z-\lambda'-\frac 12 \zeta \rho(\lambda')|^2dm(z)=\\
  \int_D|g(z)|^2|z-\lambda'|^2dm(z)+ \frac{ \rho(\lambda')^2}4
            \int_D|g(z)|^2 dm(z).     
\end{multline*}
It remains to use the relation
$$
\int_D|g(z)|^2 dm(z) \asymp \rho(\lambda')^{-2} \int_D|g(z)|^2|z-\lambda'|^2dm(z).
$$
to get \eqref{02}.

However, $\|f_1|\Lambda\|_{ \elll^2(1/\|{\bf k}_\lambda\| )}\ge CN
\|f|\Lambda\|_{ \elll^2(1/\|{\bf k}_\lambda\| )}$, and we get a contradiction
to \eqref{02} for large $N$.
\smallskip

(b) Suppose that $\dist(z,\Lambda)\ge N\rho(z)$.
Consider the function $G$ from the proof of Lemma~\rm\ref{le1a}.
Using Lemmas~\rm\ref{le0}, \rm\ref{le1a}, and \rm\ref{le2}~(a) we obtain that
for large $N$,
\begin{gather*}
\|G|\Lambda\|^2_{\elll^2(1/\|{\bf k}_\lambda\| )} = 
     \sum_{\lambda \in \Lambda} \frac{|G(\lambda)|^2}{\|{\bf k}_\lambda \|^2_\varphi}
\prec  \sum_{\lambda \in \Lambda} |G(\lambda)|^2
e^{-2\varphi(\lambda)}(\rho(\lambda))^2\\
\prec
\sum_{\lambda \in \Lambda} \int_{|\zeta - \lambda|< \rho(\lambda)}
|G(\zeta)|^2e^{-2\varphi(\zeta)}dm(\zeta)
\\
\prec 
\int_{|\zeta-z|>N\rho(z)/2} |G(w)|^2
e^{-2\varphi(w)}dm(w) \prec \frac{1}{N^2(\rho(z))^2}.
\end{gather*}
This contradicts to \eqref{03} for large $N$.
\end{proof}

\begin{proof}[Proof of Theorem~\rm\ref{te1}] 

Suppose that the system 
$\{{\Bbbk}_\lambda \}_{\lambda \in \Lambda}$  is  a Riesz basis
for $\Fp$.
 Relation \eqref{eq05} and Lemma~\rm\ref{le1a} imply 
that
$$
\Bigl\| \frac{E}{\cdot-\lambda} \Bigr\|^2_{\varphi}\asymp
|E'(\lambda)|^2\rho^2(\lambda)e^{-2\varphi(\lambda)}, \qquad \lambda \in \Lambda. 
$$
Consider the function $E/[(\cdot-\lambda)\Phi_\lambda]$. Applying  the mean value property 
we obtain
$$
\int_{\CC}\frac{|E(z)|^2}{|z-\lambda|^2}e^{-2\varphi(z)}dm(z)
\prec \rho^2(\lambda) \frac 1{\rho^4(\lambda)}
\int_{|z-\lambda|<\rho(\lambda)}
|E(z)|^2e^{-2\varphi(z)}dm(z).
$$
Take large $N$ and $w\in\mathbb C$ such that 
\begin{equation}
\label{05a}
N\ll |w|/\rho(w).
\end{equation}
By Lemma~\ref{le2}~(a), we have
\begin{gather*}
\int_{\CC}\Bigl(\sum_
{|\lambda-w|<N\rho(w), \, \lambda \in \Lambda}
\frac{1}{|z-\lambda|^2}\Bigr)|E(z)|^2e^{-2\varphi(z)}dm(z)\\
\prec
\int_{|w-z|<(N+2)\rho(w)}
\rho^{-2}(z)|E(z)|^2e^{-2\varphi(z)}dm(z),
\end{gather*}
and, as a result,
\begin{equation}
\inf_{z:|z-w|<(N+2)\rho(w)}
\Bigl[\rho^{2}(z)\Bigl(\sum_
{|\lambda-w|< N\rho(w),\, \lambda \in \Lambda}
\frac{1}{|z-\lambda|^2}\Bigr)\Bigr]\prec 1.
\label{06}
\end{equation}

Finally, by Lemma~\ref{le2}~(b),
$$
\rho^{2}(z)\Bigl(\sum_{|\lambda-w|< N\rho(w),\, \lambda \in \Lambda}
\frac{1}{|z-\lambda|^2}\Bigr)\succ
     \int_{\rho(w)<|\zeta|<N\rho(w)}
  \frac{dm(\zeta)}{|z-\zeta|^2}.
$$
For large $N$ we get a contradiction to \eqref{06}.
\end{proof}

\medskip

\noindent{\bf Remark}  {\em The above construction became possible  due to the fact  that one can choose a sufficiently large $N$ satisfying \eqref{05a}. This is not the case 
in the situations considered in Parts {\rm B}  and {\rm C}.  }
\bigskip

\section{Proofs. Part B.}
\label{sect4}

\begin{proof}[Proof of Lemma~\rm\ref{le3}] 
If $\Lambda=\{\lambda_n\}_{n\ge 0}$, 
$|\lambda_n|=\exp [(n+1)/2]$, $|z|=\exp t$, and 
$$
\frac m2-\frac 14\le t< \frac m2+\frac 14,
$$
then
\begin{gather*}
\log|E(z)|=\sum_{0\le k<m-1}\log\frac{|z|}{|\lambda_k|}
+\log\bigl|1-\frac z{\lambda_{m-1}}\bigr|+O(1)\\=
\sum_{0\le k< m}(t-(k+1)/2)+\log\dist(z,\Lambda)-t+O(1)
\\=mt-m(m+1)/4+\log\dist(z,\Lambda)-t+O(1)
\\=
t^2-3t/2+\log\dist(z,\Lambda)+O(1), \qquad |z|\to\infty.
\end{gather*}
\end{proof}
\smallskip

We  will also need the  estimate
\begin{equation}
 \label{eq06}
|E'(\lambda_n)|\asymp |\lambda_n|^{-3/2}e^{\varphi(\lambda_n)},
\end{equation}
which can be obtained similarly.
\smallskip

\begin{proof}[Proof of Lemma~\rm\ref{le4}] 
Given $\lambda\in\CC$ such that $r=\log|\lambda|\ge 1$, 
choose $n$ such that $n\le 2r<n+1$, and set $f(z)=z^n$. Then
\begin{gather*}
\|f\|^2_{\varphi}=\int_{\CC}|z|^{2n}e^{-2\varphi(z)}dm(z)
\asymp \int_0^\infty e^{2(n+1)t-2t^2}dt\\
= e^{(n+1)^2/2}\int_{-(n+1)/2}^\infty e^{-2s^2}ds
\prec e^{2(n+1)r-2r^2}.
\end{gather*}
Furthermore, $|f(\lambda)|=e^{nr}$, and we conclude that
$$
\|{\bf k}_\lambda\|_\varphi^2\succ e^{2r^2-2r}. 
$$

The opposite inequality is proved as in Lemma~\rm\ref{le1a}.
\end{proof}

\begin{proof}[Proof of Theorem~\rm\ref{te2}] 
It suffices to prove that the mapping $f \mapsto f|\Lambda$ is an isomorphism between 
$\Fp$  and $ \elll^2(1/\|{\bf k}_\lambda\|)$. 

It is straightforward that this mapping is bounded. Indeed, if $f\in \Fp$, then 
the mean value theorem and Lemma~\ref{le0} yield
\[
  |f(\lambda_k)|^2 e^{-2\varphi(\lambda_k)}|\lambda_k|^2 \prec
       \int_{|z-\lambda_k|<|\lambda_k|/10} 
              |f(z)|^2 e^{-2\varphi(z)}dm(z),
 \]
Since the discs $|z-\lambda_k|<|\lambda_k|/10$  are disjoint we obtain 
\begin{multline*}
\|f|\Lambda\|_{ \elll^2(1/\|{\bf k}_\lambda\|)}^2 \\
\prec
\sum_k   \int_{|z-\lambda_k|<|\lambda_k|/10} 
              |f(z)|^2 e^{-2\varphi(z)}dm(z) \le
 \|f\|_\varphi^2, \qquad f \in \Fp.
\end{multline*}

It is also straightforward that the mapping $f \mapsto f|\Lambda$ has zero kernel. Were this not the case we could take a non-zero $f\in \Fp$  which vanishes on $\Lambda$
and note that by Lemmas \ref{le3}  and \ref{le4}  the entire function $g=f/E$ satisfies
$|g(z)|\prec 1+|z|^{1/2}$ if  $\dist(z,\Lambda)> |z|/10$.  The latter restriction can be removed just by the maximum principle, so, by the Liouville theorem we have
$g(z)=C$ or $f(z)=CE(z)$  for some constant $C$.  Now we see that $C=0$, otherwise 
$f\not\in \Fp$, thus arriving to a contradiction.
     
It remains to prove that the mapping $f \mapsto f|\Lambda$ acts onto 
$ \elll^2(1/\|{\bf k}_\lambda\|)$, i.e. the interpolation problem
\eqref{eq03a}  has a (unique) solution for each $\{a_\lambda\}\in \elll^2(\Lambda)$.

For finite sequences  $\{a_\lambda\}$ the solution is given by the mapping
\begin{equation*}
\label{eq07}
T_\Lambda : \{a_\lambda\} \ \mapsto \ T_\Lambda\{a_\lambda\}(z)=
       \sum_{\lambda\in \Lambda} a_\lambda \|{\bf k}_\lambda \|_\varphi
             \frac {E(z)}{E'(\lambda)(z-\lambda)}.
\end{equation*}             
We will prove that 
\begin{equation}
\label{eq08}
\| T_\Lambda\{a_\lambda\} \|_\varphi \prec \|\{a_\lambda\}\|_{\elll^2(\Lambda)},
\end{equation}
and then $T_\Lambda$ extends continuously to the whole $\elll^2(\Lambda)$.

Denote
\[
E_\lambda(z)= 
    \|{\bf k}_\lambda \|_\varphi
             \frac {E(z)}{E'(\lambda)(z-\lambda)}, \qquad   \lambda \in \Lambda.
\]
Relation \eqref{eq08}  obviously follows from the inequalities 
\begin{equation}
\label{eq09}
\left \| E_\lambda  \right \|_\varphi \prec 1, \qquad 
                     \lambda \in \Lambda,
\end{equation}
and, for some $c>0$, 
\begin{equation}
\label{eq10}
\bigl| \langle E_{\lambda_m},  E_{\lambda_n} \rangle_\varphi \bigr| \prec 
                      e^{-c|n-m|}, \qquad m,n \ge 0.
\end{equation}

We are now proving these inequalities. It follows from Lemma~\rm\ref{le4} and 
\eqref{eq06} that
\[
\frac{\|{\bf k}_{\lambda_n} \|_\varphi}
     {|E'(\lambda)|}  \asymp |\lambda|^{1/2}, \qquad \lambda \in \Lambda.
\]
Together with \eqref{eq11}  this yields     
\[
|E_\lambda(z) |\asymp
        |\lambda|^{1/2} \frac{\dist(z,\Lambda)}{(1+|z|^{3/2})|z-\lambda|}
               e^{\varphi (z)},
\]
 and
 \begin{gather*}
 \|E_\lambda \|_\varphi^2 \asymp    
       \int_\CC \frac {|\lambda| \dist(z,\Lambda)^2}
               {(1+|z|^{3})|z-\lambda|^2} dm(z)  \\
= \left \{
      \int_ {|z|\leq |\lambda|/2} + 
                  \int_{|\lambda|/2<|z|<2|\lambda|}+ 
                              \int_{|z|>2|\lambda|}             
\right \} 
      \frac {|\lambda| \dist(z,\Lambda)^2}
               {(1+|z|^{3})|z-\lambda|^2} dm(z) \\
               = I_1+I_2+I_3.
\end{gather*} 

For $|z|<|\lambda|/2$  we use that $\dist(z,\Lambda) \prec |z| $, 
$|z-\lambda|\asymp |\lambda|$ to get
\[
I_1\prec \frac 1 {|\lambda|} \int_{|z|<|\lambda|/2} \frac 1 {1+|z|}  dm(z) \asymp 1.
\]

For $|\lambda|/2 < |z| < 2|\lambda|$  we use that
 $\dist(z,\Lambda)\prec |z-\lambda|$, $|z|\asymp |\lambda|$ to get
 \[
 I_2 \asymp \frac 1 {|\lambda|^2} \int_{|\lambda|/2< |z|< 2|\lambda|}  dm(z)
        \asymp 1.
 \]
 
Finally, when $|z|>2|\lambda|$  we use that $\dist(z, \Lambda)  \prec |z|$,
$|z-\lambda|\asymp |z|$ to get 
\[
I_3 \prec |\lambda| \int_{|z|>2|\lambda|} \frac{dm(z)}{1+|z|^3}  \asymp 1.
\]
 
A combination of these estimates yields \eqref{eq09}.

Furthermore, if $0\le k<m$, $\gamma=e^{1/4}$, then
\begin{gather*}
|\langle E_{\lambda_k},E_{\lambda_m}\rangle_{\Fp} |
\prec
\int_{\CC}\frac{\dist(z,\Lambda)^2}{1+|z|^3}\cdot
\frac{|\lambda_k|^{1/2}|\lambda_m|^{1/2}}
{|\lambda_k-z||\lambda_m-z|}\,dm(z)\\
=\int_{|z|<|\lambda_k|/\gamma}\ldots+
\int_{|\lambda_k|/\gamma<|z|<\gamma|\lambda_k|}\ldots+
\int_{\gamma|\lambda_k|<|z|<|\lambda_m|/\gamma}\ldots\\+
\int_{|\lambda_m|/\gamma<|z|<\gamma|\lambda_m|}\ldots+
\int_{|z|>\gamma|\lambda_m|}\ldots
\asymp \frac{|\lambda_k|^{1/2}}{|\lambda_m|^{1/2}}
\bigl(1+\log |\lambda_m/\lambda_k|^{1/2}\bigr)\\
\le ce^{-|k-m|/5},
\end{gather*}
which gives \eqref{eq10}.
\end{proof}
\bigskip

\section{Proofs. Part C.}
\label{sect5}

\begin{proof}[Proof of Lemma~\rm\ref{le5}] 
We have
$$
e^{w_n}=\int_0^\infty r^{2n}e^{-2(\log^+r)^{1+\delta}}2\pi r\,dr\asymp
\int_0^\infty e^{(2n+2)s-2s^{1+\delta}}ds.
$$
Thus, to prove the lemma, we need to describe the asymptotic
behavior of the integral
$$
\int_0^\infty e^{as-s^{1+\delta}}ds
$$ 
as  $a\to+\infty$.

This  can be done by applying standard asymptotic techniques, so we omit
calculations. 
\end{proof}

\begin{proof}[Proof of Theorem~\rm\ref{te3}] 
We need the following auxiliary statement

\begin{lem}
The numbers $r_n=\exp[(w_{n+1}-w_{n-1})/4]$ satisfy the inequality
\begin{equation}
r_n^{2n}e^{-w_n}\succ {(n+1)^2(s+1)^2}r_n^{2s}e^{-w_s} ,
\qquad n>0,\,s\not=n.
\label{07}
\end{equation}
\end{lem}
\begin{proof}
It suffices to prove that, for each $A>0$, there exists a constant $C$ such that   
\begin{multline}
\label{eq12}
\frac{n-s}2 \left [
(n+2)^{1+1/\delta} - n^{1+1/\delta} 
                \right]
       +(s+1)^{1+1/\delta} - (n+1)^{1+1/\delta} \geq  \\
 C+A|n-s|\log(n+2)+ A\log(s+2), \quad  n>0,\,s\ge 0, \, n\ne s.               
\end{multline}
Then inequality \eqref{07} will follow from Lemma ~\rm\ref{le5}.

Let $n>s$. Denote $\omega(t)=t^{1+1/\delta}$. We have 

\begin{multline*} 
 \frac{n-s}2 \left [
(n+2)^{1+1/\delta} - n^{1+1/\delta} 
                \right]
       +(s+1)^{1+1/\delta} - (n+1)^{1+1/\delta}= \\
\frac{n-s}2 \int_n^{n+2}\omega'(t)dt - \int_{s+1}^{n+1} \omega'(t)dt \succ \\
\left \{ \int_{n+1}^{n+2}\omega'(t)dt -     
     \int_{n}^{n+1}\omega'(t)dt
\right \} + 
    \left \{   \int_{n}^{n+1}\omega'(t)dt -     
      \int_{s+1}^{s+2}\omega'(t)dt
\right \}  \succ \\
n^{-1+1/\delta} + \left \{ \begin{array}{rl}
(n-s-1)n^{-1+1/\delta}  &  \text{if} \  s>n/2 \\
n^{1/\delta} & \text{if}  \ s\le n/2.
                                      \end{array}
                            \right .          
\end{multline*}
Relation \eqref{eq12}  now follows. The case $s>n$ is treated in a similar way.
\end{proof}
 
Let  now $n>0$, $\lambda_n=r_ne^{i\theta_n}$.
Then
\begin{gather*}
\Bigl\|\frac {z^ne^{-in\theta_n}}{\|z^n\|_{\varphi}} -\Bbbk_{\lambda_n}\Bigr\|_{\varphi}^2\\ =
\frac{\sum_{s\ge 0,\, s\not=n}r_n^{2s}e^{-w_s}}
{\sum_{s\ge 0}r_n^{2s}e^{-w_s}}
+\Bigl| e^{-w_n/2}-
\frac{r_n^{n}e^{-w_n}}
{(\sum_{s\ge 0}r_n^{2s}e^{-w_s})^{1/2}}
\Bigr|^2e^{w_n}\\=S_1+S_2.
\end{gather*}

By \eqref{07},
$$
S_1\le \sum_{s\ge 0,\, s\not=n}\frac{c}{(n+1)^2(s+1)^2}
\le \frac{c}{(n+1)^2}.
$$
Furthermore,
$$
S_2=\Bigl| 1-\frac{r_n^{n}e^{-w_n/2}}
{(\sum_{s\ge 0}r_n^{2s}e^{-w_s})^{1/2}}\Bigr|^2
\le 
1-\frac{r_n^{2n}e^{-w_n}}{\sum_{s\ge 0}r_n^{2s}e^{-w_s}}=S_1. 
$$
Thus,
$$
\sum_{n\ge 0}\Bigl\|\frac {z^ne^{-in\theta_n}}{\|z^n\|_{\varphi}} -\Bbbk_{\lambda_n}\Bigr\|_{\varphi}^2<\infty.
$$
By the Bari theorem (see \cite[section A.5.7.1]{NI1}), for some 
$N<\infty$, the linear operator
$U$ on $\Fp$ determined by the equalities $U(z^n)=z^n$, $0\le n<N$, 
$U(z^n/\|z^n\|_\varphi)=\Bbbk_{\lambda_n}$, $n\ge N$, extends to an isomorphism.
Therefore, the system 
\begin{equation}
\label{eq13}
\{1,z,\ldots,z^{N-1},\Bbbk_{\lambda_N},\Bbbk_{\lambda_{N+1}},\ldots\}
\end{equation}
is a Riesz basis in $\Fp$. 

Let $f\in\Fp$ be orthogonal to all $\Bbbk_{\lambda_n}$, $n\ge 0$.
We define $g\in\Fp$ by
$$
g(z)=\prod_{1\le n<N}\frac{z}{z-\lambda_n}f(z).
$$
Then $g$ is orthogonal to all the elements of the Riesz basis \eqref{eq13}
which is impossible. Therefore, the system $\{\Bbbk_{\lambda_n}\}_{n\ge 0}$
is complete in $\Fp$. Since the system $\{\Bbbk_{\lambda_n}\}_{n\ge N}$ is a Riesz basis
in a subspace of codimension $N$, another application of the Bari theorem yields that
the system $\{\Bbbk_{\lambda_n}\}_{n\ge 0}$ is a Riesz basis in $\Fp$.
\end{proof}

\bigskip
\bigskip
\bigskip

\noindent \textsc{Alexander Borichev, Centre de Math\'ematiques et Informatique,
Universit\'e d'Aix-Marseille I,
39 rue Fr\'ed\'eric} 

\noindent \textsc{Joliot-Curie,
13453 Marseille,
France}

\noindent\textsl{E-mail}: \texttt{borichev@cmi.univ-mrs.fr}
\medskip

\noindent \textsc{Yurii Lyubarskii, Department of Mathematical
Sciences, Norwegian University of Science and Technology,
N-7491 Trondheim, Norway}

\noindent\textsl{E-mail}: \texttt{yura@math.ntnu.no}

\end{document}